\documentstyle[12pt]{article}
\newcommand{\braid}[3]{{#1}$\lower4pt\hbox{$\oo\atop\raise4pt
           \hbox{$\scriptscriptstyle {#3} $}$}${#2}}
\newcommand{\twist}[3]{{#1}${\,\scriptscriptstyle {#3}}\atop\raise9pt
           \hbox{$\scriptstyle\oo$} ${#2}}
\newcommand{\ve}{\varepsilon}
\newcommand{\be}{\begin{eqnarray}}
\newcommand{\ee}{\end{eqnarray}}
\newcommand{\n}{\nonumber }
\newcommand{\oo}{\otimes}
\newcommand{\bt}{\beta}

\newcommand{\si}{\sigma}
\newcommand{\gm}{\gamma}

\newcommand{\hH}{\hat H}
\newcommand{\hX}{\hat X}
\newcommand{\hE}{\hat E}
\newcommand{\tX}{\tilde X}
\newcommand{\al}{\alpha}

\def\C{{\mathchoice {\setbox0=\hbox{$\displaystyle\rm C$}\hbox{\hbox
to0pt{\kern0.4\wd0\vrule height0.9\ht0\hss}\box0}}
{\setbox0=\hbox{$\textstyle\rm C$}\hbox{\hbox
to0pt{\kern0.4\wd0\vrule height0.9\ht0\hss}\box0}}
{\setbox0=\hbox{$\scriptstyle\rm C$}\hbox{\hbox
to0pt{\kern0.4\wd0\vrule height0.9\ht0\hss}\box0}}
{\setbox0=\hbox{$\scriptscriptstyle\rm C$}\hbox{\hbox
to0pt{\kern0.4\wd0\vrule height0.9\ht0\hss}\box0}}}}

\begin{document}
\begin{titlepage}

\begin{center}
{\Large \bf Twisting cocycles in fundamental representaion
and triangular bicrossproduct Hopf algebras
 }
\end{center}
\vspace{0.5in}
\begin{center}
{\Large  A.I. Mudrov
 }
\end{center}
\begin{center}
\vspace{0.2in}
 Department of Theoretical Physics,
Institute of Physics, St.Petersburg State University, Ulyanovskaya 1,
Stary Petergof, St. Petersburg, 198904, Russia\\
(aimudrov@DG2062.spb.edu)
\end{center}
\vspace {0.5in}
\begin{center}
{\bf Abstract}
\end{center}
{\footnotesize
We find the general solution to the twisting equation in the tensor bialgebra
$T({\bf R})$ of an associative unital ring ${\bf R}$ viewed as
that of fundamental representation for a universal enveloping Lie algebra and
its quantum deformations. We suggest a procedure of constructing twisting
cocycles
belonging to a given quasitriangular subbialgebra ${\cal H}\subset T({\bf R})$.
This algorithm generalizes Reshetikhin's approach, which involves cocycles
fulfilling the Yang-Baxter equation. Within this framework we study a class of
quantized inhomogeneous Lie algebras related to associative rings in a certain
way, for which we build twisting cocycles and universal $R$-matrices.
Our approach is a generalization of the methods developed for
the case of commutative rings in our recent work
including such well-known examples as Jordanian
quantization of the Borel subalgebra of $sl(2)$  and the null-plane quantized
Poincar\'e algebra by Ballesteros {\em at al}.
We reveal the role of special group cohomologies in this process
and establish the bicrossproduct structure of the examples studied.
}
\vspace{0.5in}
\begin{center}
{1998}
\end{center}

\end{titlepage}

\section{Introduction}
Quantum deformations of Lie groups and algebras are at present a
subject of intensive studies from the viewpoints of collecting facts and
crystallizing mathematical concepts as well as
of searching for new physical applications. Among the
established notions of the quantum group theory one should
mention Drinfeld's twisting \cite{D1,D2} and Majid's
bicrossproduct and doublecrossproduct constructions \cite{Mj4}.
Twisting, realizing a specific equivalence
between two Hopf algebras plays an important role
for the geometrical and physical reasons because
it controls deformation not only of the symmetry algebra
of a manifold but of its whole geometry coherently
Therefore a classification  of quantum deformations of, say,
a universal enveloping Lie algebra ought to provide the answer
about twist-equivalence between its different types.
Majid's doublecrossproduct construction has close connection
with twisting and the quantum double in particular \cite{D3,RSTS,Mj1}.
As for the bicrossproduct, its relation to
quasitriangularity and twisting
is not so well understood, despite of numerous examples
including quasitriangular Hopf algebras.
The most significant step in that direction was made in
Ref. \cite{Mj2}, where the double of the algebra
$\C({\bf M})\bowtie \C{\bf G}$, built on a matched pair of groups
${\bf M}$ and ${\bf G}$ was shown to be a bicrossproduct itself.
As examples of  bicrossproduct we would like to mention
the $\kappa$-deformation of the Poincar\'e algebra
\cite{LRTN,Mj3}, the canonical example of the Jordanian
quantization of the Borel subalgebra of  $sl(2)$ \cite{Mj4},
and the null-plane quantized Poincar\'e algebra \cite{BHOS,AHO}.
The last two are the results of twisting of classical
universal enveloping algebras and are therefore triangular.
As we have shown in Ref. \cite{M1}, they  are associated with
commutative rings\footnote{Throughout the paper by ring we mean
a finite-dimensional algebra over a field ${\bf K}$ so as to reserve
the word "algebra" for Hopf one.}, one-dimensional in the first
case and that
spanned by matrices
$$
\al_1=\left(
\begin{array}{lll}
0 & 0 & 1\\
0 & 0 & 0\\
0 & 0 & 0\\
\end{array}
\right),\quad
\al_2=\left(
\begin{array}{lll}
0 & 0 & 0\\
0 & 0 & 1\\
0 & 0 & 0\\
\end{array}
\right),\quad
\al_3=\left(
\begin{array}{lll}
1 & 0 & 0\\
0 & 1 & 0\\
0 & 0 & 1\\
\end{array}
\right)
$$
in the second. Along this line we have found a generalization
of the examples mentioned for an arbitrary commutative ring.
The algebras  studied in Ref. \cite{M1} were twisted  classical universal
enveloping algebras and also bicrossproduct Hopf algebras, although not
considered in that context. In the present paper we formulate a generalization
of our approach to an arbitrary associative ring providing a new class
of quasitriangular bicrossproduct Hopf algebras.
The Hopf operations are explicitly written out in terms
of generators and the twisting 2-cocycles and universal
$R$-matrices are presented. The class under investigation
arises as an example of twisting which may be
regarded as a generalized Reshetikhin's procedure \cite{R}
involving a solution to the Yang-Baxter equation as a
twisting cocycle. We come to this generalization analyzing
solutions to the twisting equation in the tensor bialgebra
$T({\bf R})$ of an associative unital ring ${\bf R}$, considered as
that of the fundamental representation of a given Hopf algebra
${\cal H}$ (we call a homomorphism ${\cal H} \to {\bf R}$
fundamental if its lifting to whole $T({\bf R})$, which
always existits, is non-degenerate).

The paper is organized as follows. Section II is auxiliary and
contains general description of the tensor bialgebra $T({\bf R})$ structure,
its subbialgebras and homomorphisms. This part may be considered
as the "non-coordinate" formulation of the Faddeev-Reshetikhin-Takhtajan
method \cite{RTF} suitable for an arbitrary associative unital ring.
Since rings with identities have exact matrix realizations, e. g.
by the regular representations on themselves, such a reformulation does not
supply
with particular new information compared with the traditional matrix
approach. Nevertheless, it provides certain technical convenience, so we find
it possible to present this formulation here. Section III is devoted to
solving the twisting equation in $T({\bf R})$. Therein we develop a procedure
of constructing "universal" cocycles starting from elements of
${\bf R}^{\oo 2}$
obeying certain conditions. This algorithm is illustrated in Section IV
on inhomogeneous Lie algebras related in a sense to associative
rings. We build the deformed coproduct, quantum commutation relations,
find twisting cocycles and universal  $R$-matrices.
In Section V the connection between the investigated algebras
and the bicrossproduct construction is established.

\section{Bialgebra $T({\bf R})$}
 To perform algebraic manipulations it is convenient to formulate
 the algorithm by Faddeev, Reshetikhin, and Takhtajan  \cite{RTF} in
 the part of
 constructing quantum algebra Fun$_q({\bf R})$ of functions on matrix
 rings for an arbitrary associative ring ${\bf R}$. We need some information
 concerning its structure and the structure of homomorphisms
 from a Hopf
 algebra or, more generally, bialgebra ${\cal H}$ into
 Fun$_q({\bf R})$. Let $\mu$ be the multiplication in ${\bf R}$.
 We choose a basis $(x^{\alpha}) \subset {\bf R}$ such that
 $ x^0 $ is the identity of ${\bf R}$.
 The dual basis in  ${\bf R}^*$ will be marked with subscripts.
 Denote $F({\bf R})$ the algebra over a field ${\bf K}$ freely
 generated by 1 and  $(x_\alpha)\in {\bf R}^*$. Introduce the
 coproduct and the counit defining them  on the generators as
$$ \Delta (x_\alpha) = \mu_{\alpha}^{\rho \sigma} x_{\rho}\oo x_{\sigma},\quad
\Delta (1) = 1\oo 1, \quad $$
$$\epsilon(1) = 1,\quad \epsilon(x_0) = 1,\quad  \epsilon(x_i) = 0,\quad
i\not=0$$
and extending over whole $F({\bf R})$ homomorphically.
The dual bialgebra $T({\bf R})=F^*({\bf R})$ appears to be a direct
sum of its ideals
$\sum^{\infty}_{n=0} {\bf R}^{\oo n} $, where
$ {\bf R}^0 $ coincides the field ${\bf K}$ of scalars.
The multiplication in $T({\bf R})$ is characterized by the property
${\bf R}^{\oo n} {\bf R}^{\oo m} = 0$ for $m \not = n$.
The identity of $T({\bf R})$ is expanded as the sum $\sum^{\infty}_{n=0}e^n$
of idempotents, where $e^0$ is the unity of ${\bf K}$, and $e^n$ with $n>0$
are the those of ${\bf R}^{\oo n}$.
Multiplying by $e^n$ carries out the projection homomorphism
\begin{displaymath}
\begin{array}{ccc}
      \pi^m\colon T({\bf R}) &\to & {\bf R}^{\oo m},
\end{array}
\end{displaymath}
and for $n=0$ this is just the bialgebra counit.
The coproduct in $T({\bf R})$ is determined by the product of its dual
algebra $F({\bf R})$ and on the basis elements is defined by the formula
\be
  \Delta (x^{i_1..i_n}) &=& e^0\oo x^{i_1..i_n}+...+x^{i_1..i_k}
  \oo x^{i_{k+1}..i_n}+...+x^{i_1..i_n}\oo e^0,
  \quad x^{i_1..i_n}\in {\bf R}^{\oo n}.\n
\ee
It follows from here that the composite mapping
\begin{displaymath}
\begin{array}{ccccccc}
   {\bf R}^{\oo(i+j)} &\to & T({\bf R}) &\stackrel{\Delta}{\to} &
   T({\bf R})\oo T({\bf R}) &\to &{\bf R}^{\oo i}\oo {\bf R}^{\oo j},
\end{array}
\end{displaymath}
where the left arrow means the injection and the right one is
the projection homomorphism, turns out to be a ring
isomorphism. Quotient of $F({\bf R})$ by the ideal ${J}$ generated by
quadratic relations of the form
$ x_{\alpha} x_{\beta} = B^{\rho \sigma}_{\beta \alpha} x_{\rho} x_{\sigma}$
inherits the coproduct if and only if the subspace
${\bf S} \in {\bf R}^{\oo 2} $
of functionals annihilating these relations is a subalgebra in
${\bf R}^{\oo 2}$. In particular, such  a subalgebra can be
determined as the set of solutions to the equation
\be
 \ R z R^{-1}&=& \tau z , \quad z\in {\bf R}^{\oo 2}, \n
\ee
($\tau$ permutes the factors in ${\bf R}^{\oo 2}$),
and then the relations in the dual algebra will look as
$$ R_{\gm \nu}\mu_{\alpha}^{\gm \rho}
               \mu_{\beta}^{\nu \sigma} x_{\rho} x_{\sigma} =
                x_{\sigma} x_{\rho} \mu_{\alpha}^{\rho\gm}
 \mu_{\beta}^{\sigma\nu}R_{\gm \nu}.$$
The bialgebra ${\cal U}$  dual to the factor-bialgebra
$ {\cal A} \equiv F({\bf R}) / {J} $ is decomposed into the direct sum of
its ideals
$ \sum^{\infty}_{n=0} \pi^n ({\cal U}) $, where  $\pi^1 ({\cal U})$ is
isomorphic to the ring ${\bf R }$ itself and each addend $\pi^n({\cal U})$
at $n>1$ is a subring in $ {\bf R}^{\oo n}$
characterized by the condition
$z_{...\al\bt ...}=B^{\rho\si}_{\bt\al}z_{... \rho\si ...}$
такая for $ z \in  \pi^n({\cal U})$.
Evidently, $\pi^n({\cal U})$ is the intersection of
all possible subalgebras
${\bf R}^{\oo i} \oo \pi^2({\cal U})\oo {\bf R}^{\oo j}$,
such that $i+j+2=n$.

Let us describe the structure of a homomorphism $\phi$ of an arbitrary
bialgebra ${\cal H}$ to ${\cal U}$. The composition of $\phi$ with the
projector $\pi^n$ is an algebraic mapping.
We set $\rho\equiv\pi^1 \circ\phi$ and $\rho^n\equiv\pi^n \circ\phi$ and
also introduce the notations
$$\Delta^1\equiv id\colon {\cal H}\to {\cal H},\quad
\Delta^2\equiv \Delta\colon {\cal H}\to {\cal H} \oo {\cal H},\quad
\Delta^3\equiv (\Delta\oo id) \Delta\colon {\cal H}
\to {\cal H}\oo {\cal H}\oo {\cal H}, ...$$
Let $h\in {\cal H}$ and $x_{\vec \al} \in {\cal A}$, where  $\vec \al$
is a multiindex of length $n$.
>From the chain of equalities
  \be
  \langle \rho^n (h),x_{\vec \al}\rangle &=&
  \langle \pi^n \circ\phi (h),x_{\vec \al}\rangle =
  \langle \phi (h),x_{\vec \al}\rangle =
  \langle \Delta^n\circ\phi (h),x_{\al_1}\oo...\oo x_{\al_n} \rangle =
  \n\\&&
  \langle \phi^{\oo n}\circ \Delta^n (h),x_{\al_1}\oo...\oo x_{\al_n}\rangle =
  \n\\&&
  \langle (\pi^1\circ\phi)^{\oo n}\circ \Delta^n(h),
  x_{\al_1}\oo...\oo x_{\al_n}\rangle =
  \langle \rho^{\oo n}\circ \Delta^n (h),x_{\al_1}\oo...\oo x_{\al_n}\rangle
  \n
  \ee
we find that conditions
\begin{enumerate}
\item  $\rho^n = \rho^{\oo n} \circ \Delta^n$,
\label{it1}
\item $\rho^2({\cal H}) \subset {\bf S}=\pi^2({\cal U})$,
\label{it2}
\end{enumerate}
are satisfied, and these completely specify $\phi$.
Namely, if an algebraic mapping   $\rho\colon{\cal H} \to {\bf R}$
fulfills condition \ref{it2}, there exists the unique bialgebra
homomorphism of $\phi \colon {\cal H} \to {\cal U}$
such that  $\pi^n\circ\phi=\rho^n$.
Indeed, mapping $\rho^{\oo n} \circ\Delta^n$ respects the algebraic
structure for every $n$. Obviously,
$\rho^{\oo (n+2)} \circ\Delta^{n+2} ({\cal H}) \subset
  {\bf R}^{\oo i} \oo {\bf S}\oo {\bf R}^{\oo(n-i)}$ for every $i\leq n$,
because of the coassociativity of the coproduct. We define
$\phi$ via the formula
$\phi(h) = e^0 \ve(h)+\sum_{n> 0}\rho^{\oo n}\circ\Delta^n(h)$ and verify
that it is also a coalgebra homomorphism
(its algebraic property is evident).
For two multiindices ${\vec \al}$ and ${\vec \bt}$ of lengths $m$ and $n$
we have
$$
  \langle (\phi \oo \phi) \Delta (h),x_{\vec \al}\oo x_{\vec \bt}\rangle =
  \langle (\rho\oo...\oo\rho) \Delta^{m+n} (h),
   x_{\al_1}\oo...\oo x_{\al_m}\oo x_{\bt_1}\oo...\oo x_{\bt_n}\rangle =
$$
$$
\>
   =
  \langle \phi (h),x_{\vec \al}x_{\vec \bt}\rangle =
  \langle \Delta\circ\phi (h),x_{\vec \al}\oo x_{\vec \bt}\rangle \n
$$
  that proves what required.
  Conditions \ref{it1} and \ref{it2} are convenient practical criteria
  for checking out the homomorphic properties of
  mappings from ${\cal H}$ to ${\cal U}$.

\section{Twisting equation in $T({\bf R})$}
Algebras we are interested in, namely, universal enveloping Lie algebras and
their quantum deformations appear to be embedded into
$T({\bf R})$ associated with their fundamental representation ring ${\bf R}$.
Since twisting of a subalgebra induces that of the whole
algebra, the problem can be put forward of describing all
solutions to the twisting equations in $T({\bf R})$
with the hope of further selecting among them those belonging to
required subalgebras. We can point out the following advantages
of such an approach. As far as the composition of two twistings is concerned,
thus we avoid inconvenience  of dealing with a bialgebra different from
original after performing the first deformation.
Another remarkable feature of such a description is a possibility
to reduce the elaborate task of constructing the universal
twistor to the much easier problem of solving
a system of equations, rather non-linear, in a finite-dimensional
ring. Our study resulted in finding the general solution
to the twisting equation in $T({\bf R})$. We have managed to formulate
conditions more general than those employed in Reshetikhin's approach
\cite{R} which ensure that the twisting cocycle built on its image in
${\bf R}^{\oo 2}$
would lie in the required subalgebra.

Consider the twisting equation
\be
 (\Delta\oo id)(\Phi)\Phi_{12}= (id\oo\Delta)(\Phi)\Phi_{23},
\label{TE}
\ee
in the tensor cube of a bialgebra ${\cal H}$
where the subscripts determine the embeddings
${\cal H}^{\oo 2}\to {\cal H}^{\oo 3}$.
An invertible solution to this equation can be normalized
in such a way that
\be
(\ve\oo id)(\Phi)=(id\oo\ve)(\Phi) =1.
\label{Norm}
\ee
Such a solution takes part in transforming
 bialgebra ${\cal H}$ into a new one with the same multiplication
and the
coproduct $\tilde\Delta(h)=\Phi^{-1}\Delta(h)\Phi$,
$h\in {\cal H}$. Other objects, e.g. the counit, antipode, universal
${\cal R}$-matrix, if any, are connected with the old ones via
the well-known formulas which can be found in \cite{D2}.
We are going to prove the following assertion.
\newtheorem{theor}{Theorem}
\begin{theor}
For every set of invertible elements $\Phi^{1,k}\in {\bf R}\oo {\bf R}^k$,
there exists the unique solution  $\Phi \in T({\bf R})\oo T({\bf R})$ to the
twisting equation, such that  $(\pi^1\oo \pi^k)(\Phi)=\Phi^{1,k}$.
\end{theor}
Having applied the projector $\pi^m\oo\pi^n\oo\pi^k$ to the both
sides of equality (\ref{TE}), we come to the equation in
${\bf R}^{\oo (m+n+k)}$:
\be
\Phi^{m+n,k}\Phi_b^{m,n}= \Phi^{m,n+k}\Phi^{n,k}_e,
\label{TE1}
\ee
where $\Phi^{m,n}$ is the image of the twisting cocycle
$(\pi^m\oo\pi^n)(\Phi)$. Letters $b$ and $e$ indicate that
the elements $\Phi_b^{m,n}$ and $\Phi_e^{n,k}$ are embedded into
${\bf R}^{\oo (m+n+k)}$ from the "beginning" and from the "end", respectively:
${\Phi_b^{m,n}} \in {\bf R}^{\oo (m+n)}\oo e^k \subset {\bf R}^{\oo (m+n+k)}$,
${\Phi_e^{n,k}} \in e^m\oo {\bf R}^{\oo (n+k)} \subset {\bf R}^{\oo (m+n+k)}$.
Because of (\ref{Norm}), for every
$k\geq 0$ we have $\Phi^{0,k} = e^0\oo e^k$ and $\Phi^{k,0} =  e^k \oo e^0$.
Suppose now that the  elements $\Phi^{1,k}$ are known. Then,
using  equation (\ref{TE1}), $\Phi^{m,k}$ can be defined
recursively for all $m$ and $k$ greater than 1:
\be
 \Phi^{m+1,k}\equiv \Phi^{m,1+k}\Phi^{1,k}_e\bar{\Phi}_b^{m,1}
 \label{build}
\ee
(the bar stands for the inverse). This implies the uniqueness of
the solution. Obviously, equation (\ref{TE1}) is true when one of
the numbers $m$, $n$, and $k$ are equal to zero. By construction,
it is fulfilled for $n=1$ and all $m$ as well. So we must show
that equation (\ref{TE1}) is satisfied for arbitrary $m$, $n$, and $k$.
Assume the required property proved for all $m$ and $n$, which sum
is less than $N>2$. Then for $m+n=N$ and $k>1$ we have
$$ \Phi^{m+n,k}= \Phi^{m+n-1,1+k}\Phi^{1,k}_e\bar{\Phi}_b^{m+n-1,1}.$$
Within the assumption made, we decompose the first factor on the
right-hand side of this equality, according to (\ref{build}),
and rewrite (\ref{TE1}) in the equivalent form
$$\Phi^{m,n+k}\Phi^{n-1,1+k}_e\bar{\Phi}_b^{m,n-1}\Phi^{1,k}_e
\bar{\Phi}_b^{m+n-1,1}=\Phi^{m,n+k}\Phi^{n,k}_e\bar{\Phi}_b^{m,n}.$$
Dividing both sides by the first factor
we come to condition
$$\Phi^{n-1,1+k}_e\bar{\Phi}_b^{m,n-1}\Phi^{1,k}_e
\bar{\Phi}_b^{m+n-1,1}=\Phi^{n,k}_e\bar{\Phi}_b^{m,n}.$$
Again, decomposing the first factor on the right-hand side
according to the recursion assumption  we find
$$\Phi^{n-1,1+k}_e\bar{\Phi}_b^{m,n-1}\Phi^{1,k}_e\bar{\Phi}_b^{m+n-1,1}=
\Phi^{n-1,1+k}_e\Phi^{1,k}_e\bar{\Phi}_{\{m+1\}}^{n-1,1}\bar{\Phi}_b^{m,n},$$
where the subscript in $\bar{\Phi}_{\{m+1\}}^{n-1,1}$  means that
it is embedded into ${\bf R}^{\oo (m+n+k)}$ beginning from $\{m+1\}$-th
place. It is important that the factor
$\bar{\Phi}_b^{m,n-1}$ on the left-hand side can be permuted with
$\Phi^{1,k}_e$. Division by $\Phi^{n-1,1+k}_e\Phi^{1,k}_e$ yields
$$\bar{\Phi}_b^{m,n-1}\bar{\Phi}_b^{m+n-1,1}=
\bar{\Phi}_{\{m+1\}}^{n-1,1}\bar{\Phi}_b^{m,n}$$
or
$${\Phi}_b^{m+n-1,1}{\Phi}_b^{m,n-1}=
{\Phi}_b^{m,n}{\Phi}_{\{m+1\}}^{n-1,1}.$$
This is exactly twisting equation (\ref{TE1}) for
$k\to 1$ and $n\to n-1$. According to the induction principle we
consider the theorem proved.

Thus, the  family of solutions to equation (\ref{TE}) in the bialgebra
$T({\bf R})$
turns out to be very large: it is parameterized by an arbitrary set
of invertible elements  $\Phi^{1,k}$, $k>0$. On the other hand,
we are interested only in those $\Phi$ which belong to the subalgebra
${\cal U}\oo {\cal U}$. We cannot propose the general method to build
such solutions. However, we can point out an algorithm which
can help to solve the problem at least for quasitriangular ${\cal U}$
and which can be interpreted as a generalization of Reshetikhin's
approach. Starting form $\Phi^{1,1}$ as known, set
$\Phi^{1,k}\equiv \Phi_{12}\Phi_{13}\ldots\Phi_{1 (k+1)}$.
Here $\Phi_{1i}$ are the images of  $\Phi^{1,1}$ via the corresponding
embeddings ${\bf R}^{\oo 2} \to {\bf R}\oo {\bf R}^{\oo k}$.
For $\Phi^{1,k}$ to belong to ${\cal U}\oo {\cal U}$, it is necessary
and sufficient to require
\be
R_{23}\Phi_{12}\Phi_{13}=\Phi_{13}\Phi_{12}R_{23}.
\label{Split2}
\ee
Further, element $\Phi^{2,1}$ lies in ${\cal U}\oo {\cal U}$ if and only if
$R_{12}\Phi^{2,1}=(\tau\oo id)(\Phi^{2,1})R_{12}$.
Having expressed $\Phi^{2,1}$ through $\Phi^{1,1}$ we come
to equation
$R_{12}\Phi_{12}\Phi_{13}\Phi_{23}\bar{\Phi}_{12}=
       \Phi_{21}\Phi_{23}\Phi_{13}\bar{\Phi}_{21} R_{12}$
or, the matrix
$\tilde{R}=\bar\Phi_{21}R\Phi$ introduced, to equation
\be
\tilde{R}_{12}\Phi_{13}\Phi_{23}=\Phi_{23}\Phi_{13}\tilde{R}_{12}.
\label{Split1}
\ee
Verification of the condition
 $R_{12}\Phi^{2,k}=(\tau\oo id)(\Phi^{2,k})R_{12}$ for $k>1$
boils down to   equality
$\tilde{R}_{12}\Phi_{13}\ldots\Phi_{1\> k+2}
               \Phi_{23}\ldots\Phi_{2\> k+2}=
               \Phi_{23}\ldots\Phi_{2\> k+2}
               \Phi_{13}\ldots\Phi_{1\> k+2}\tilde{R}_{12}$,
which is, as can be easily seen, follows from (\ref{Split1}).
Now, with further use of identity (\ref{TE1}) one can
see that $\Phi^{m,n}$ belongs to ${\cal U}\oo {\cal U}$
for every $m$ and $n$, indeed.

\section{Twisting of inhomogeneous Lie algebras}
One can notice that the example of solution to the twisting equation
built at the end of the previous section satisfies the identities
\be
(id \oo\Delta)(\Phi)=\Phi_{12}\Phi_{13},
\label{RRP1}
\ee
\be
(\tilde\Delta\oo id)(\Phi)=\Phi_{13}\Phi_{23}.
\label{RRP2}
\ee
where $\tilde\Delta$ is the twisted coproduct:
$\tilde\Delta(h)=\bar\Phi\Delta(h)\Phi$. Reshetikhin's conditions
are obtained from here if $\Phi$ solves the Yang-Baxter equation
and, besides, $(\Delta\oo id)(\Phi)=\Phi_{23}\Phi_{13}$.
This generalization of Reshetikhin's twisting is non-trivial,
the non-standard quantization of the Borel subalgebra of $sl(2)$
taken into account. The twisting cocycle for
$U(b(2))$ has the form $\exp({X\oo H})$, where $H$ is the primitive
element in $U(b(2))$ and $X$ is the primitive element in
$U_h(b(2))$ \cite{Ogiv}. This example was generalized in
Ref. \cite{M1} for an arbitrary commutative ring
which in the case of $U(b(2))$  coincides with the field of
scalars.
In the present section we shall formulate the analogous
generalization for an arbitrary associative ring ${\bf L}$, not
necessarily commutative.
Let the multiplication in ${\bf L}$ be defined by the structure
constants $B_{\mu\nu}^\si$. Consider a Lie algebra
built on $H_\mu \in {\bf L}$ and $X^\nu \in {\bf L}^*$ subjected
to the commutation relations
$$[H_\mu,H_\nu]=(B_{\mu\nu}^\si-B_{\nu\mu}^\si) H_\si,\quad
[H_\mu,X^\nu]=-B_{\mu\si}^\nu X^\si.$$
The subalgebra ${\bf L}^*$ generated by $X^\nu$ is assumed to be
commutative. It is easy to see that the element
$H_\nu\oo X^\nu- X^\nu\oo H_\nu$ satisfies the classical Yang-Baxter
equation. The subspace ${\bf L}^*$ is a right module over ring
${\bf L}$. Let us affiliate the identity to ${\bf L}\triangleright {\bf L}^*$
and denote the resulting ring ${\bf R}$. The product in  ${\bf R}$
is evaluated according to the rules
$$\hH_\mu \hH_\nu=B_{\mu\nu}^\si \hH_\si,\quad
\hX^\nu\hH_\mu =B_{\mu\si}^\nu \hX^\si,\quad
\hX^\mu \hX^\nu=0,\quad
 \hH_\mu \hX^\nu=0,$$
plus evident expressions involving identity $\hat E$.
Starting from this multiplication, one can see that the element
\be
\tilde R=\hat E\oo \hat E + \hat H_\nu\oo \hat X^\nu -
\hat X^\nu\oo \hat H_\nu
\label{QYBE}
\ee
is a solution to the quantum Yang-Baxter equation, and the
element
\be
\Phi^{1,1}= \hE\oo \hE - \hX^\nu\oo \hH_\nu
\label{Cocycle}
\ee
obeys (\ref{Split2}) and (\ref{Split1}) where
$R$ should be set to $\hE\oo \hE$. Correspondence
$1 \to \hat E$, $X \to \hat X$, $H \to \hat H$  is extended to
a homomorphism (non-degenerate) of the universal enveloping algebra
$U({\bf R})$ into $T({\bf R})$.
\begin{theor}
   Twisting cocycle $\Phi$ expanding $\Phi^{1,1}$ (Eq. (\ref{Cocycle}))
   by modified Reshetikhin's procedure belongs to $U({\bf L}^*)\oo U({\bf L})$
   and has the formд $\exp(-\tX^\nu \oo H_\nu)$, where
   $\tX^\nu$ are expressed by series in  $X^\nu$.
\end{theor}
Indeed, $\Phi^{1,k} \in U({\bf L}^*)\oo U({\bf L})$  by construction.
>From  defining formula (\ref{build}), recursively using the
facts that $Span(\hE,\hH_\mu)$ and $Span(\hX^\mu)$ form a subring
and an ideal in ${\bf R}$, respectively, we get the first assertion of
the theorem. Now the announced form of the twisting cocycle
follows from (\ref{RRP1}).

Unemployed so far identity (\ref{RRP2}) enables us to determine
the twisted coproduct on the generators $\tX^\nu$. It turns out that
$\tilde\Delta(\tX^\nu)=D^\nu(1\oo\tX,\tX\oo 1)$, where
$D(a,b)$ is the Campbell-Hausdorf series: $e^{D(a,b)}=e^a e^b$ for
arbitrary $a$ and $b$ from the Lie algebra of the ring ${\bf L}$.
Thus, the commutative algebra $U_q({\bf L}^*)$ is isomorphic to the function
algebra on the group $\exp({\bf L})$ taken with the opposite coproduct.
We have yet to evaluate the twisted coproduct on $H_\nu$ and
to determine commutation relations $[H_\nu,\tX^\mu]=f(\tX)^\mu_\nu$.
What can be said about functions $f(x)^\mu_\nu$ is that they
are subject to the "boundary" conditions
свойством
\be f(0)^\mu_\nu=0, \quad
\partial_\si f(0)^\mu_\nu= - B^\mu_{\nu\si}.
\label{boundcond}
\ee
It is accounted for the following. First of all, knowing the image
$\Phi^{1,1}$ of the cocycle $\Phi$ in ${\bf R}^{\oo 2}$ we conclude
that $\frac{\partial\tX^\mu}{\partial X^\nu}|_{x=0}=\delta^\mu_\nu$.
Now the required properties of $f(x)^\mu_\nu$ are conditioned by
the homomorphism from $U({\bf R})$ to  ${\bf R}$ and the strong
nilpotence of $\hX^\mu$.

Having introduced matrices of the left and right regular representations
$L(X)^\mu_\nu \equiv B^\mu_{\si\nu}X^\si$,
$R(X)^\mu_\nu \equiv B^\mu_{\nu\si}X^\si$, from the definition of
$\tilde\Delta$ we find
\begin{eqnarray}
\tilde\Delta(H_\mu)&=&\exp(\tX\oo H)(H_\mu\oo 1 + 1\oo H_\mu)
\exp(-\tX\oo H) \n\\
   &=&H_\mu\oo 1 - \Biggl(\frac{e^{L(\tX)-R(\tX)}-1}
   {\scriptstyle L(\tX)-R(\tX)}f(\tX)\Biggr)^\nu_\mu\oo H_\nu\n\\
   &+&(e^{L(\tX)-R(\tX)})^\nu_\mu \oo H_\nu
   \n\\
   &=&H_\mu\oo 1 + g(\tX)^\nu_\mu \oo H_\nu .
   \label{Copr}
\end{eqnarray}
The coassociativity requirement imposed, formula (\ref{Copr}) implies
$g(D(1\oo\tX,\tX\oo 1))=g(1\oo\tX)g(\tX\oo 1)$, and this
with necessity  entails $g(\tX)=e^{A(\tX)}$, where $A$ represents
a left action of ${\bf L}$ on itself.
Resolving $g(\tX)$ with respect to $f(\tX)$ and using conditions
(\ref{boundcond}) we come finally to
$A(\tX)^\mu_\nu=L(\tX)^\mu_\nu=B^\mu_{\si\nu}X^\si$.
The resulting formulas describing the Hopf structure of the twisted
algebra $U_q({\bf R})$ read
\be
\tilde\Delta(\tX^\mu)&=&D(1\oo \tX,\tX\oo 1)           ,\n\\
\tilde\Delta(H_\mu)&=&H_\mu\oo 1+(e^{L(\tX)})^\nu_\mu \oo H_\nu  ,\n\\[5pt]
[H_\mu, H_\nu]&=&(B_{\mu\nu}^\si-B_{\nu\mu}^\si) H_\si,    \n\\[5pt]
[H_\nu,\tX^\mu]&=&
\Biggl(\frac{\scriptstyle L(\tX)-R(\tX)}{e^{L(\tX)-R(\tX)}-1}
\Bigl(e^{L(\tX)-R(\tX)}-e^{L(\tX)}\Bigr)\Biggr)^\mu_\nu.\n\\
\label{relHX0}
\ee
The antipode is easily found from the coproduct:
$$S(\tX^\mu)=-\tX^\mu,\quad S(H_\mu)=-(e^{-L(\tX)})^\nu_\mu H_\nu.$$
The expressions obtained generalize formulas deduced
for commutative ring in Ref. \cite{M1}. In the latter case
the value of the commutator $[H_\nu,\tX^\mu]$ is simplified because of
${L(\tX)=R(\tX)}$ and turns into $(1-e^{L(\tX)})$,
as for the well-known example of the Borel subalgebra of  $sl(2)$.
The universal $R$-matrix is expressed through the twisting cocycle
by the standard formula ${\cal R}=\Phi^{-1}_{21}\Phi$ \cite{D1}:
\be
{\cal R}&=&\exp( H_\nu\oo\tX^\nu )\exp(-\tX^\nu \oo H_\nu)
\label{URM}
\ee
and has the form familiar from the theories of the Jordanian
quantization of $sl(2)$ \cite{Ogiv} and the null-plane quantized
Poincar\'e algebra \cite{M1,BHP}.
Thus we obtain a closed and complete description of the
deformed algebra $U_q({\bf R})$, although it would be desirable
to find the relation between
$\tX^\mu$ and the classical generators $X^\mu$.
To this end, let us calculate the twisted coproduct on elements
$X^\mu$:
\be
 \tilde\Delta(X^\mu)
 &=&\exp(\tX\oo H)(1\oo X^\mu + X^\mu\oo 1 )\exp(-\tX\oo H) \n\\
 &=&X^\mu\oo 1+ (e^{-L(\tX)})^\mu_\nu\oo X^\nu,\n
\ee
that results in the functional equation
\be
\varphi(\Delta(\tX))&=&\varphi(\tX')+ e^{-L(\tX')}\varphi(\tX''),
\label{gr_coc}
\ee
where the primes distinguish the tensor components, and  $\varphi$
is the transformation connecting the quantum and classical generators:
$X^\mu=\varphi^\mu(\tX)$. This equation is well known from the theory
of the group cohomologies \cite{G} and its solution is
$$\varphi(\tX)=\frac{e^{-L(\tX)}-1}{\scriptstyle-L(\tX)}\tX.$$
This formula solves  the problem of proceeding to the classical basis of
$U_q({\bf R})$.

The analysis of the solution found allows us to perform
a further generalization of the examples considered above
in the following direction. Let ${\bf G}$ be a Lie
group and ${\bf L}$ its Lie algebra with the basis elements $H_\mu$.
Assume a left action $H_\mu\triangleright H_\nu=B_{\mu\nu}^\si H_\si$
of ${\bf L}$ on itself, which is, as a rule, does not
coincide with the adjoint representation. Let function
$\varphi\colon {\bf G}\to {\bf L}$ be a group 1-cocycle,
that is $\varphi(ba)=\varphi(a)+ a^{-1}\triangleright\varphi(b)$,
$a,b \in {\bf G}$.  It can be viewed
as a mapping defined in some neighborhood of the origin in
${\bf L}$. We suppose $\varphi$ to be invertible and denote its
inverse $\psi$. By the left conjugate action
$H_\mu\triangleright X^\nu=-B_{\mu\si}^\nu X^\si$
on the dual space we build the semidirect sum
${\bf L}\triangleright {\bf L}^*$, where ${\bf L}^*$
is considered as an Abelian subalgebra.
\begin{theor}
   Element $\Phi=\exp(-\psi^\nu(X) \oo H_\nu)$ is a twisting
   cocycle for the universal enveloping Hopf algebra
   $U({\bf L}\triangleright {\bf L}^*)$.
\end{theor}
Notice that $\exp(-\psi^\nu(X) \oo H_\nu)$ satisfies identity
(\ref{RRP1}). Hence, in order to prove the theorem
the second identity (\ref{RRP2}) should be stated.
Making use of the fact that $\tilde\Delta$ is an algebraic
mapping for arbitrary $\Phi$, we evaluate $\tilde\Delta(\tX)$
on the elements $\tX^\mu=\psi^\mu(X)$ and come to the equation
(\ref{gr_coc}), where operator $L(\tX)$ is defined via
tensor $B$ as above. Because of the invertibility
of  $\varphi$, this implies $\Delta(\tX^\mu) = D^\mu(1\oo\tX,\tX\oo 1)$.
Then identity (\ref{RRP2}) is obeyed as well.

\section{Bicrossproduct structure}
To exhibit the bicrossproduct structure of the examples considered
let us use the FRT method and recover the quantum groups
by the solution to the quantum Yang-Baxter equation
(\ref{QYBE}). Modulo the order of the factors, they are isomorphic to
the quantum algebras, and the isomorphism is realized via the
universal $R$-matrix (\ref{URM}).

In terms of the basis  $(e, x_\mu,h^\nu)$ dual to the basis
$(\hE, \hH_\mu,\hX^\nu)$ of the ring ${\bf R}$ the coproduct,
according to the scheme rendered in Section II,
has the form
$$\Delta(h^\si)=h^\si\oo e + e \oo h^\si + B^\si_{\mu\nu} h^\mu\oo h^\nu,
  \quad
 \Delta(x_\si)=x_\si\oo e + e \oo x_\si+ B_{\mu\si}^\nu  x_\nu\oo h^\mu,$$
$$\Delta(e)=e\oo e.$$
The counit is determined by the rule $\ve(e)=1$, $\ve(x_\mu)=0$,
$\ve(h^\mu)=0$. Imposing RTT-type relations with the matrix $R$
given by (\ref{QYBE}), we come to the following permutaion rules
\be
[x_\mu,x_\nu]&=& (B^\si_{\nu\mu} -
B^\si_{\mu\nu}) x_\si e,\label{relXX} \\[6pt]
[x_\mu,h^\nu]&=& B^\nu_{\mu\si}h^\si e + B^\nu_{\si\al}
B^\si_{\bt\mu}h^\al h^\bt.\label{relXH}
\ee
Other commutation relations are trivial and, in particular,
the element $e$ belongs to the center of the algebra.
Note, that the ideal $(e-1)$ is a Hopf one and set $e=1$.
Introduce quantities $\eta^\mu$ starting from equality
$h^\nu\oo \hH_\mu = e^{\eta^\mu\oo \hH_\nu}-1\oo \hE$.
In terms of new generators ($\eta^\mu,x_\mu$) the coproduct
turns out to be opposite to that of the algebra $\tilde U({\bf L})$,
and that is seen through the substitution
$\eta^\mu \to  \tX^\mu$, $ x_\mu\to -H_\nu $.
Commutation relations (\ref{relXX}) are thus recovered exactly.
It has yet to be shown that relations (\ref{relXH})
goes over into the last expression in (\ref{relHX0}). This is guaranteed by
the uniqueness of the value of the commutator
(\ref{relXH}) as a function in $h^\mu$ compatible with
the given coproduct and fulfilling
$\frac{\partial[x_\mu,h^\nu]}{\partial h^\si}|_{h=0}=B^\nu_{\mu\si}$.
This boundary condition is determined by the
homomorphism from the quantum group into the ring ${\bf L}$
via the "square" matrix $R$ (\ref{QYBE}).

According to \cite{Mj4}, bicrossproduct ${\cal A}\bowtie {\cal B}$
of Hopf algebras  ${\cal A}$ and ${\cal B}$ is defined
via a left action $\triangleright$  of ${\cal A}$ on
${\cal B}$ and a right coaction $\beta$ of ${\cal B}$ on ${\cal A}$.
The latter is a mapping from ${\cal A}$ into the tensor product
${\cal A} \oo {\cal B}$. The conjugate mapping to $\beta$
realizes a right action ${\cal A}^* \triangleleft {\cal B}^*$.
These operations are subjected to the  set of consistency conditions
\cite{Mj4}.
Multiplication and comultiplication on ${\cal A}\bowtie {\cal B}$
are evaluated via
$$(a\oo h)(b\oo g)= a(h_{(1)}\triangleright b)\oo h_{(2)}g,
$$
$$\Delta (a\oo h)= (a_{(1)}\oo h_{(1)}^{\>(\bar1)}) \oo
(a_{(2)}h_{(1)}^{\>(\bar2)}\oo h_{(2)}),$$
where  $a,b \in {\cal A},\quad h,g \in {\cal B}$,
$\Delta(h)=h_{(1)}\oo h_{(2)}$,
and $\bt(h)=h^{(\bar1)} \oo h^{(\bar2)}$.
Turning to $U_q({\bf R})$ constructed in the previous section we see that
algebraically it has the same structure as
$ U^*({\bf L})\triangleleft U({\bf L})$,
whereas its dual quantum group has the form
$U^*({\bf L})\triangleright U({\bf L})$.
In terms of generators $h^\mu$ and $x_\nu$ the left action is
$x_\mu\triangleright h^\nu=B^\nu_{\mu\si}h^\si + B^\nu_{\si\al}
B^\si_{\bt\mu}h^\al h^\bt$,
and the right coaction $\beta$ can be computed from the coproduct formula
$$id\oo\bt(h)\oo id = (id\oo id\oo id\oo\ve)\circ \Delta(1\oo h).$$
This yields
$$ \bt(x_\si)=x_\si\oo 1 + B_{\mu\si}^\nu  x_\nu\oo h^\mu$$
for the generators $x^\mu$ thus stating the bicrossproduct structure of
$U^*_q({\bf R})$ and  $U_q({\bf R})$.
\section{Discussions}
The problem of explicitly calculating twisting 2-cocycles for
Hopf algebras is a non-trivial  one even if their existence
is {\em a priory} known. Difficulties arise already in the simplest
case of classical universal enveloping Lie algebras, despite of
advanced Drinfeld's theory on quantizing triangular Lie bialgebras
\cite{D1}.
This explains why examples of explicitly given twisting 2-cocycles
are in relatively short supply. So, it seems quite natural to reduce
the problem to studying "quadratic" 2-cocycles which are images
of "universal" ones in fundamental representations, provided
there exists some "fusion" procedure to expand those
matrix solutions over representations of higher spins.
There are two algorithms of this kind \cite{H,H9609029,JC}, both based on
factorisation properties of twisting elements  \cite{RSTS,R,Mj1,M2}
$(\Delta \oo id) (\Phi) = \Phi_{13}\Phi_{23}$ or
$(\Delta \oo id) (\Phi) = \Phi_{23}\Phi_{13}$ (and appropriate
identities involving $id \oo \Delta$).
Depending on the order of the factors on the right hand side,
additional requirements like $\Phi_{12}\Phi_{23}=\Phi_{23}\Phi_{12}$
or the Yang-Baxter equation are imposed on $\Phi$. Although most of explicitly
known universal twisting 2-cocycles are due to these two options,
it is clear that they cannot cover all possibilities. The idea of
proceeding to fundamental representation in a finite-dimensional
ring ${\bf R}$ in studying twist-equivalences among various quantizations
seems
yet more fruitful because the associated bialgebra $T({\bf R})$
plays the role of a container, in a general situation, for
all the deformations of a Hopf algebra. Thus every
twisting cocycle of  a subbialgebra undergoing deformation
remains so for whole $T({\bf R})$. This makes it reasonable
to consider twisting equation in $T({\bf R})$ and then
try to select solutions belonging to the  given
subbialgebra. The first part of this program has been
completely solved in the present paper, while for
the second we have suggested a new kind of fusion procedure
which appears to be close to Reshetikhin's twisting.
The novelty is that the conditions imposed on $\Phi$
employ both twisted and untwisted coproducts. To demonstrate effectiveness
of the proposed scheme we have considered a class of inhomogeneous universal
enveloping Lie algebras related to associative rings in a special way  and
quantized them along that line. The technique
can be viewed as a generalization of the theory developed
in our previous work \cite{M1} dealing exclusively with commutative rings,
which was motivated by the Jordanian quantization of $sl(2)$
and the null-plane quantized Poincar\'e algebra. We have also
exposed the bicrossproduct structure of the objects investigated
thus providing new examples of quasitriangular bicrossproduct
Hopf algebras. A remarkable fact is that the class considered
may be treated directly with the use of special Lie group
cohomologies which take part in building twisting elements.
Relevance of Hopf algebra cohomologies to
the twisting procedure and bicrossproduct construction was
already pointed out in Ref. \cite{Mj4}, so the present study
gives a new insight to their role in the theory.
It is interesting to generalize cohomological methods
applied here for  classical universal enveloping algebras
to Hopf algebras of more general nature.

\vspace{0.3cm}
\noindent
{\Large\bf  Acknowledgement}
\vspace{0.3cm}

\noindent
We are grateful to P. P. Kulish and V. D. Lyakhovsky for
helpful and stimulating discussions.

\end{document}